\documentclass[12pt]{article}
\usepackage{geometry}                
\usepackage[all]{xy}
\usepackage{graphicx}
\usepackage{amssymb}
\usepackage{epstopdf}
\begin{document}
 \newcommand{\bfz}{\textbf{0}}
  \newcommand{\bfu}{\textbf{1}}
  \newcommand{\mR}{\mathcal{R}}
   \newcommand{\mI}{\mathcal{I}}
    \newcommand{\I}{\mathcal{I}}
 \newtheorem{definition}{Definition}
\newtheorem{lemma}{Lemma}
\newtheorem{theorem}{Theorem}
\newtheorem{remark}{Remark}
\newtheorem{example}{Example}

\author{Vittorio Cafagna and Gianluca Caterina}
  \title{MV-algebras and measure: some examples}
    \maketitle

\abstract{These notes explore a class of examples of MV-algebras. Our  point of view is that the closed unit interval $\mathcal I=[0,1]$, the prototype of MV-algebras, is also a mathematical archetype, as the only (up to diffeomorphisms) connected 1-manifold with non-connected boundary. Possibly on a same level as the other archetypes: $\mathbb R$, the only non compact 1-manifold, $S^1$, the only compact 1-manifold and $\mathbb R^+$, the only non-compact 1-manifold with a connected boundary. As such, $\mathcal I$ has an ubiquitous role in many different parts of mathematics: homotopy theory, Morse theory of compact manifolds, deformation theory of complex varieties, moduli spaces, flow boxes in dynamical systems, not to mention the obvious cases of probability theory and  fuzzy theory. Nonetheless, it does not seem that the algebraic structure of $\mathcal I$  play any role in any of the fore-mentioned theories (with the  obvious exception of fuzzy theory). $\mathcal I$ is bounded, so to speak,  to the passive role of a parameter space. We propose in this article a  definition of a MV-algebra structure on a class of subsets of some probability spaces and we work-out some examples. Our intention is to convey, by mean of the  simplest possible examples, the idea that the topology of the geometric object under consideration might be reflected in the MV-algebraic structure. We also discuss an example of Chang \cite{ch} and discuss similarities and differences with the proposed class of examples. This workwas part of a broader project jointly with Professor Vittorio Cafagna which spans a few years between 2000 and 2005. Professor Cafagna passed away unexpectedly in 2007. His intellectual breadth and inspiring passion for mathematics, however, is still very well alive.}\footnote{The support of Professor Antonio Di Nola is gratefully acknowledged. Thanks to Professors J. Golan, and G. L. Litvinov for their precious comments regarding the algebraic part of this work, and to Professor Y.Katsov  for having sent me some of his papers on tensor product of semimodules.
}

\section{Introduction}

MV-algebras have been introduced by Chang \cite{ch} in the attempt to generalize some algebraic techniques from Boolean logic to multi-valued logics, from which the term MV derives. The first and, in some sense, most fundamental example of MV-algebra is the unit interval $\mathcal{I}=[0,1]$ endowed with a binary operation called \L ukasiewicz t-norm, also called {\sl truncated sum} since it captures the idea of threshold, intrinsic to interval-like structures. Truncated sum is a particular instance of a more general class of binary operations on $\mathcal{I}$, to which sometimes we refer as {\sl t-norms}, which stays for {\sl triangular norms}. The problem which has motivated the examples presented in this article has been the search for mathematical structures whose relation to $\mathcal{I}$ make them good candidates for being intrinsically associated with MV-algebras. Simple measure spaces seem to fit well in such a perspective; however, the extent of the generality up to we can extent our considerations is not clear yet, and at the end of this article we discuss a very simple example which induces what we called a `generalized' MV-algebra since the last axioms (the commutativity of the join), does not hold. As byproduct of our searching for these kind of structures we present a class of t-norms which appears to be retaining  some geometric properties of the object under consideration. 

 \section{Basic definitions}

An MV-algebra is a commutative monoid 
 \[(A,\oplus,\neg,0,1)\]
 where $\oplus$ is the structure operation, $0$ is the identity, and 
 \begin{enumerate}
 \item $\neg(\neg(x))=x$
 \item $1\oplus x=1$
 \item $x\oplus\neg 0=\neg 0$
 \item $\neg(\neg x\oplus y)\oplus y=\neg(\neg y\oplus\neg x)\oplus x$
 \end{enumerate}

\begin{example}

\begin{enumerate}
\item The archetypal model for MV-algebras is the {\sl \L ukasiewicz t-norm}
$(\mathcal{I},\oplus,\neg,0,1)$ where
 \[x\oplus y=min(x+y,1)\]  
 and
  \[\neg x=1-x\]
  \item (Chang,\cite{ch}) Define the following families of formal symbols:
   \[I_0=\{0,a,2a,3a,\dots \}\] \[I_1=\{1,\overline{a},\overline{2a},\overline{3a},\dots\}\]  
 where $\overline{ka}=1-ka$. Letting $+$ be the ordinary sum between integers, define the following binary operation $\oplus$ on $\hat{I}=I_0\cup I_1$:
  \begin{enumerate}
  \item $na\oplus ma=(n+m)a$
  \item $\overline{na}\oplus \overline{ma}=1$
  \item $ma\oplus \overline{na}=\left\{\begin{array}{c}1\ \ m\leq n \\
  \\ \overline{(m-n)a}\ m>n\end{array}\right.$
  \item $\overline{ma}\oplus na=\left\{\begin{array}{c}1\ \ m\leq n \\ \\ \overline{(m-n)a}\ m>n\end{array}\right.$
  \end{enumerate}
    Define an involution:
  \begin{enumerate}
\item $\neg ka=\overline{ka}$
\item $\neg(\overline{ka})=ka$
 \end{enumerate}
 \end{enumerate}
 Chang \cite{ch} has proved that:
 \begin{theorem}
$(\hat{I},\oplus, \neg,0,1)$ is an MV-algebra. 
\end{theorem}
\begin{flushright}
\end{flushright}
 \end{example}
Ex. 2 is very interesting in our perspective, since it somewhat abstracts the concept of `boundary' of an MV-algebra: the entire structure $\hat{I}$, as set, is defined by two special points (which are named 0 and 1 but do not have to be confused with 0 and 1 as end points of $\mathcal{I}$) and two sets $I_0$ and $I_1$ associated with them, with the property that $I_0$ is the bijective image of $I_1$ under the involution (and viceversa). MV-algebras which we are going to discuss contain strong analogies with $(\hat{I},\oplus, \neg,0,1).$
 \section{Intervals}
 Chang completeness theorem assigns to the MV algebra  $(\mathcal{I},\oplus,\neg,0,1,)$ a special role in the set of all MV-algebras, in the same way as Stone representation theorem does to  the Boolean algebra $\{0,1\}$ with respect to the set of all Boolean algebras: if an equation holds in $\mathcal{I}$, then it holds in every MV-algebra. As we have mentioned, {\textit t-norms} are a generalization of the operation of truncated sum (they are actually generalized conjunction connectives over $\mathcal{I}$). Triangular norms play a crucial role in several fields of mathematics and artificial intelligence. For an exhaustive overview on t-norms we refer to \cite{pa}. The following definition is taken from Golan \cite{go1}:
\begin{definition}
 A triangular norm ({\sl t-norm}) on $\mathcal{I}$ is a binary operation, say $+,$ on $\mathcal{I}$ satisfying the following two conditions:
\begin{enumerate}
\item ($\mathcal{I}$,+) is a commutative monoide with identity element 1.
\item $a\leq b\Rightarrow a + c\Rightarrow b + c$ for all $c\in\mathcal{I}$
\end{enumerate}
\end{definition}

The `leitmotif' of this article is the search of examples suggesting that the \L ukasiewicz t-norm can be seen as the trace of certain set-theoretic operations on some measure spaces, therefore we start with the most fundamental one. 

\begin{definition}
Let $\mathcal{A}$ the $\sigma-algebra$ over $\mathcal{I}$ generated by the closed intervals and let 
\[\mathcal{I}_0=\{[0,a]s.t.0\leq a \leq 1\}\subseteq\mathcal{A}\]
\[\mathcal{I}_1=\{[a,1]s.t.0\leq a \leq 1\}\subseteq\mathcal{A}\]
two parametrized families of subsets of $\mathcal{A}$.

\end{definition}
Let us also introduce the following symbols:
\begin{definition}
For all $a,b\in\mathcal{I}$ we set 
\[ a^0=[0,a]\in\mathcal{I}_0\]
\[a^1=[1-a,1]\in\mathcal{I}_1\]
\end{definition}
\begin{example}
The interval [0.2,1] is represented in this notation by $0.8^1$, while [0,0.6] by $0.6^0$.
\end{example}
We denote by $0$ the degenerate interval [0,0], by 1 the interval $[0,1]$; we also we denote the complement of a subset of $B\subset\mathcal{A}$ is denoted by $\overline{B}$. 
\begin{definition} 
Let $\mu$ be the Lebesgue measure on $\mathcal{A}$; we define the involution
\[i:\mathcal{I}_0\longrightarrow \mathcal{I}_0\]
\[ a^0\mapsto [\mu(\overline{a^0})]^0\]
\end{definition}
\begin{lemma}
 $i(a^0)=(\neg a)^0$ where $\neg a=1-a$.
\end{lemma}

$Proof$. $i(a^0)=[\mu(\overline  a^0)]^0=\mu[(1-a)^1]^0=(1-a)^0$.

\begin{definition}
For any $a\in \mathcal{I}$ let 
\[j:\mathcal{I}_0\longrightarrow\mathcal{I}_1\]
\[a^0\mapsto a^1\]
\end{definition}
The following lemma is a direct consequence of the definition:
\begin{lemma}
$j$ is a preserving-measure bijection, such that
$j^{-1}=j$, $j(\mathcal{I}_0)=\mathcal{I}_1$ and  $j(\mathcal{I}_1)=\mathcal{I}_0$.
\end{lemma}
$Proof.$\ $\mu[j(a^0)]=\mu(a^1)=\mu([1-a,1])=a=\mu(a^0)$. We also have that we can define 
\[j^{-1}:\mathcal{I}_1\longrightarrow\mathcal{I}_0\]
\[ a^1\mapsto a^0\]
Then we have 
\[j^{-1}[j(a^0)]=j^{-1}(a^1)=a^0\]
\begin{flushright}
$\square$
\end{flushright}

\subsection{Structure}
Let us define the following operations:

\begin{equation}
\oplus:\mathcal{I}_0\times\mathcal{I}_0\longrightarrow\mathcal{I}_0\ \ \  \ \ \ 
a^0\oplus b^0=[\mu(a^0\cup b^1)]^0
\end{equation}
\begin{equation}
 \odot:\mathcal{I}_0\times\mathcal{I}_0\longrightarrow\mathcal{I}_0\ \ \ \ \ \ 
a^0\odot b^0=[\mu(a^0\cap b^1)]^0
\end{equation}
\begin{theorem} 

Let us set $\overline{b}=1-b$; the following hold:
For any $a^0,b^0\in\mathcal{I}_0$ we have 
\begin{enumerate}
\item $a^0\cap b^1\neq\emptyset\Leftrightarrow a^0\cup b^1=\mathcal{I}$
\item $a^0\cap b^1\neq\emptyset \Rightarrow\mu(a^0\cap b^1)=(a-\overline{b})$
\item $a^0\oplus b^0=1^0\Leftrightarrow a+b\geq1 $
 \item $a^0\oplus b^0=(a+b)^0\Leftrightarrow a+b\leq1$
 \end{enumerate}
\end{theorem}

$Proof.$ \begin{enumerate}
\item$a^0\cap b^1\neq\emptyset\Leftrightarrow a\geq \overline{b}\Leftrightarrow  a^0\cup b^1=\mathcal{I}$
\item $a^0\cap b^1\neq\emptyset\Rightarrow a\geq \overline{b}\Rightarrow a^0\cap b^1=[a,\overline{b}]\Rightarrow \mu(a^0\cap b^1)=(a-\overline{b})$

\item If $a^0\oplus b^0=1$, then 
\[ [\mu(a^0\cup b^1)]^0=1\Leftrightarrow \mu(a^0\cup b^1)=1\Leftrightarrow (a\geq \overline{b})\Leftrightarrow a+b\geq 1\]

\item $a+ b\leq 1\Leftrightarrow a\leq \overline{b}\Leftrightarrow a^0\cap b^1=\emptyset\Leftrightarrow a^0\oplus b^0=\mu(a^0\cup b^1)^0=[\mu(a^0)+\mu(b^1)]^0=(a+b)^0$
\end{enumerate}
\begin{flushright}
$\square$
\end{flushright}
Let us also notice that $I_0$ is parameterized by $\mathcal{I}$ with fiber consisting of one singe element, that is 
\[\phi:\mathcal{I}\longrightarrow I_0\]
\[a\mapsto a^0\]
is a bijection; therefore the trace on $\mathcal{I}$ of the operation we have defined is the  \L ukasiewicz t-norm, so we can conclude that  $(\mathcal{I}_0,\oplus,\neg,0,1)$ is an MV-algebra.

The reader may object that the above observation makes our construction trivial, since the effect of the operation we have defined, as an arithmetic on $\mathcal{I}$, coincide with the  \L ukasiewicz t-norm. 

However, the point that seems interesting to us, is the fact that such an operation has been obtained without referring directly to the algebra on [0,1], whereas it has been induced by an operation on sets in which the measure seems to force the arithmetic. 

This is the reason why we think it is worthwhile to show that $(\mathcal{I}_0,\oplus,\neg,0,1)$ is an MV-algebra by using only set-theoretic considerations as follows: 
 \begin{theorem}
 $(\mathcal{I}_0,\oplus,\neg,0,1)$ is an MV-algebra.
 \end{theorem}
 In what follows we use the fact that $\mu$, being a finite measure, satisfies $\mu(A\cup B)=\mu(A)+\mu(B)+\mu(A\cap B)$, and the fact that, by definition, $\mu(a^0)=\mu(a^1)$.
 
$Proof.$\ It is immediate to check that $\neg 0=1$ and that, for any $a^0\in\mathcal{I}_0$, $0\oplus a^0=$ and $1\oplus a^0=1$.\\ 
Let us now notice that
\[a^0\cup b^1=\mathcal{I}\Leftrightarrow a^0\cap b^1\neq\emptyset\]
and
\[a^0\cap b^1\neq\emptyset\Rightarrow \mu(a^0\cap b^1)=(a-\overline{b})\]

We have that
\[a^0\oplus b^0=\{\mu[a^0\cup b^1]\}^0\]
since $\mu$ is a measure 
\[=\{[(\mu(a^0)+\mu (b^1))-\mu(a^0\cap b^1))]\}^0=\]
\[=\left (a+b- \left\{\begin{array}{c} (a-\overline{b}) \ \ \  if\  a^0\cap b^1\neq\emptyset \\0\ \ if\  a^0\cap b^1=\emptyset \end{array}\right.\right)^0\]
On the other hand we have that 
\[b^0\oplus a^0=\{\mu[b^0\cup j(a^0)]\}^0\]
since $\mu$ is a measure 
\[=\{[(\mu(b^0)+\mu (a^1))-\mu(b^0\cap a^1))]\}^0=\]
\[=\left( b+a- \left\{\begin{array}{c} (b-\overline{a}) \ \ \  if\  b^0\cap a^1\neq\emptyset \\0\ \ if\  b^0\cap a^1=\emptyset \end{array}\right.\right)^0\]
Since 
\[ a^0\cap b^1\neq\emptyset\Rightarrow b^0\cap a^1\neq\emptyset\]
and 
\[ b-\overline{a}=a-\overline{b}\]
then the identity
\[a^0\oplus b^0=b^0\oplus a^0\]
follows, hence $\oplus$ is commutative.\\

Let us set 
\[(A)=(a^0\oplus b^0)\oplus c^0=\{\mu[a^0\cup b^1)]\}^0\oplus c^0=\]
\[               = \mu( \{\mu[a^0\cup b^1)]\}^0)^0\cup c^1))^0   \]
and 
\[(B)=a^0\oplus (b^0\oplus c^0)=\{\mu[b^0\cup c^1)]\}^0\oplus a^0=\]
\[                \mu( \{\mu[b^0\cup c^1)]\}^0)^0\cup a^1))^0   \]

 Suppose $a^0\cap b^1=\emptyset$\\
 In this case we have:\\
 \[(A)= \mu( \{\mu(a^0)+ \mu(b^1))\}^0)^0\cup c^1))^0 \]
If we have 
\[ (a^0\oplus b^0)\cap c^1=\emptyset \]
then 
\[(A)=[((\mu(a^0)+\mu(b^1))+\mu(c^1)]^0;\]

But this is equivalent to say that $a^0\cap c^1=\emptyset$ and $b^0\cap c^1=\emptyset$; by applying these identities to $(B)$ we can conclude that 
\[(A)=(B)\]
The analysis of the remaining cases is similar, from which associativity follows.

Since $(\mathcal{I}_0,\oplus,\neg,0,1)$ is totally ordered we can  conclude that $(\mathcal{I}_0,\oplus,\neg,0,1)$ is an $MV-$ algebra.
\begin{flushright}
$\square$
\end{flushright}

Needless to say, the above discussion can be reproduced to show that:
\begin{theorem}
 $(\mathcal{I}_0,\odot,\neg,1,0)$ is an MV-algebra.
 \end{theorem}

 \subsection{Analogies with  Ex.2}

In order to highlight some analogies between the previous example and the one proposed by Chang, we first need to rearrange the notation:
\[b^1=[b,1]\]

We can define an operation in $\mathcal{I}_1$ completely similar to the one introduced for $\mathcal{I}_0$, say $\oplus_1$, as follows:
\begin{enumerate}
\item $j_1:\mathcal{I}_1\longrightarrow\mathcal{I}_0$
\[b^1\mapsto (1-b)^0\]
 \item $\neg_1(b^1)=(1-b)^1$
\item $a^1\oplus b^1= \mu[a^1\cup j_1(b^1)]$
\end{enumerate}
 Let also denote by 0 the degenerate interval [1,1] and by 1 the entire $\mathcal{I}$. It is immediate to verify that $(\mathcal{I}_1,\oplus_1,\neg_1,0,1)$ is a MV-algebra. 
 
 Let us not we define:  
\[\hat{\mathcal{I}}=\mathcal{I}_0\cup{\mathcal{I}_1}\]
 and the following binary operation $\star$ on it:
 \begin{enumerate}
 \item $a^0\star b^0=a^0\oplus b^0$
 \item $a^0\star b^1=a^1\oplus' b^1$
 \item $a^1\star b^0=a^1\oplus' b^1$
  \item $a^1\star b^1=1$
  \item $\neg a^0=a^1$
  \item $\neg a^1=a^0$
 \end{enumerate}

Let us know notice that, for any $m,n\in\mathcal{I}$,  by the above definition we have that:
\begin{enumerate}
\item $m^0\star n^0=(m+n)^0$ if $m+n\leq 1$
  \item $a^1\star b^1=1$
\item $n^0\star m^1=\left\{\begin{array}{c}1\ \ \ \  m\leq n \\
  \\ (1-(m-n))^1\ \ m>n \end{array}\right.$
  \item $m^1\star n^0=\left\{\begin{array}{c}1\ \ \ \  m\leq n \\
  \\ (1-(m-n))^1\ \ m>n \end{array}\right.$
  \end{enumerate}
 
 The formal analogy with the the Ex.2 is evident, since 1), 2), 3) and 4) corresponds to $a), b), c)$ and $d)$ of Ex.2; however, the sets that we have called $I_0$ and $I_1$ are countably infinite, in particular they do not have an upper bound, whereas $\mathcal{I}_0$ and $\mathcal{I}_1$ are uncountably infinite and bounded by $\mathcal{I}$. We would have liked to go beyond the formal analogies, and show that $\mathcal{I}$, with the operations we have just defined, is an MV-algebra, but we have encountered some problems in the proof of that, so that we leave it as an open question.

 \subsection{Rectangles}
In this section we will introduce the first non trivial generalization of the ideas of the previous section. In order to better explicit them, we first start by looking at the two dimensional analogue of the intervals:
\begin{definition}
Let $Q=[0,1]\times[0,1]$, let $\mathcal{R}$ the sigma algebra of closed balls of $Q$;  the $lower\ rectangles$ is the family of subset parameterized by $\mathcal{I}$
\[R^0(\lambda)=\{(x,y)\in Q\ |\ y\leq \lambda\}\]
and the $upper\ rectangles$ the family of subsets parameterized by $\mathcal{I}$
\[R^1(\lambda)=\{(x,y)\in Q\ |\ y\geq 1-\lambda\};\]
we also set 

\[\mathcal{R}^0=\{R^0(\lambda)\ 0\leq \lambda\leq 1\}\]
\[\mathcal{R}^1=\{R^1(\lambda)\ 0\leq \lambda\leq 1\}\]
\end{definition}
We write 0 for the degenerate rectangle $R^0(0)$ and 1 for the whole $Q=R^0(1)$.
Following the strategy used for the intervals, we can extend the operation $\oplus$  and the involution to rectangles  in a natural way:
\begin{definition}
 if $\mu$ is Lebesgue measure on $\mathcal{R}$ we can define the involution
\[i:\mathcal{R}_0\longrightarrow \mathcal{R}_0\]
\[ R^0(a)\mapsto [\mu(\overline {R^0})]^0=R^0(\neg a)\]
where $\neg a=1-a$
and the map
\[j:\mathcal{R}_0\longrightarrow\mathcal{R}_1\]
\[R^0(a)\mapsto R^1(a)\]
\end{definition}

\begin{lemma} $j^{-1}=j$, $j(\mathcal{R}_0)=\mathcal{R}_1$ and  $j(\mathcal{I}_1)=\mathcal{I}_0$.\
\end{lemma}
$Proof$. The same as in lemma 5.

\begin{definition}
We can define the binary operation
\[\oplus:\mathcal{R}_0\times\mathcal{R}_0\longrightarrow\mathcal{R}_0\]
\[R^0(a)\oplus R^0(b)=[\mu(R^0(a)\cup R^1(b))]^0\]
 and its dual
 \[\odot:\mathcal{R}_0\times\mathcal{R}_0\longrightarrow\mathcal{R}_0\]
\[R^0(a)\odot R^0(b)=[\mu(R^0(a)\cap R^1 (b))]^0\]
\end{definition}

The same argument we have used for the intervals can be adapted to show that:
\begin{theorem}
$(\mathcal{R}^0,\oplus, \neg, 0, 1)$ and $(\mathcal{R}^0,\odot, \neg, 1, 0)$ are MV-algebras.
\end{theorem}

Let us now introduce a `perturbation' in our object: let 
$Q=[0,1]\times[0,1]$, and let $Q_1$ be the interior of the square inscribed in it with `center' $(0.5,0.5)$ and edges of length $k$; let us define the following quantities:
\begin{definition} 
\begin{enumerate}
\item $\phi:\mathcal{I}\longrightarrow \mathcal{I}$
\[\lambda\mapsto \mu(\{(x,y)\in Q\ |\ y\leq \lambda\}-Q_1)\]

\item $T^0(\phi(\lambda))=\{(x,y)\in Q\ |\ y\leq \lambda\}-Q_1$
\item $T^1(\phi(\lambda))=\{(x,y)\in Q\ |\ y\geq 1-\lambda\}-Q_1$

\item $\mathcal{T}^0=\{T^0(\phi(\lambda))\ \  0\leq \lambda\leq 1\}$
\item $\mathcal{T}^1=\{T^1(\phi(\lambda))\ \  0\leq \lambda\leq 1\}$
\end{enumerate}
\end{definition}
 \begin{definition}
Let $\neg_k(x)=(1-k^2)-x$, then we define
\begin{enumerate}
\item $\neg T^0(\phi(\lambda))=T^0(\neg_k(\phi(\lambda)))$
\item $j[T^0(\phi(\lambda))]=T^1(\phi(\lambda))$

\item $T^0(\phi(a))\oplus T^0(\phi(b))=\mu[T^0(\phi(a))\cup T^1(\phi(b))]^0$
\end{enumerate}
\end{definition}

\begin{theorem}
 $(\mathcal{T}^0,\oplus,\neg,0,1)$ is an MV-algebra.
\end{theorem}

$Proof.$ First notice that the elements $T(t)\in\mathcal{T}^0$ are parameterized by their measure  $t\in[0,1-k^2]$ and the operation we have defined can be written as 
\[T^0(\phi(a))\oplus T^0(\phi(b))=T^0(\phi(a)\oplus_k \phi(b))\]
where $a\oplus_k b=min(1-k^2,a+b)$, which, along with the involution $\neg$, 0 and $K=1-k^2$, defines the MV algebra $([0,1-k^2],\oplus_k,\neg, 0,K)$. Since $\mathcal{T}^0$ is totally ordered by its measure-parameter, and since the operation we have defined on $\mathcal{T}^0$ coincides with $\oplus_k$, we can conclude that $(\mathcal{T}^0,\oplus,\neg, 0,1)$ is an MV algebra.

\section{Some Examples}

There exists a large literature  on the subject of t-norms \cite{pa},\cite{je1},\cite {je2}, approached from different point of view, \cite{ge}, and representing a t-norm as a surface in $\mathbb{R}^3$ is a tool often used to present a dynamical picture of their analysis.

 We can look at the trace on $\mathcal{I}$ of such an operation and check that is:
\[a\oplus b=\phi^{-1}(\mu((T^0(\phi(a)\oplus T^0(\phi(b)))))\]
 
 In what follows we give a picture of the t-norm induced by the simple geometric construction we have introduced.

We have noticed that some geometric characteristic of the `hole'  survive in the associated graph; while we do not have developed a solid theoretical background to give structural results about this sort of phenomena, and we have intended to present this work as a nutshell of examples to be better understood and thoroughly justified, we think that this is an indication of the role that MV-algebras could play in topology.

\subsection{The Square}
In the following example we start with the unit square with a smaller square-hole with side of length $k$ and centered, at the point $(0.5,0.5)$. Using Mathematica, we can give a dynamic picture of the deformation: as we could expect, for small values of $k$ the induced {\sl t-norm} is very close, to the \L ukasiewicz t-norm, whereas $k$ increases we can notice how the shape of the hole is somehow reflected in such a deformation.

\subsection{The Disk}
We can replace the square hole with a one of circular shape. Once again, the deformation that we obtain by increasing the length of the radius of the hole reflects the shape of the hole.

\section{Generalized MV algebras}

D. Mundici \cite{mu1} gives a functorial characterization of MV-algebras by showing the equivalence between the category  of MV algebras and $l-$groups with a strong unit. This results enforces the interpretation of MV-algebras as an abstraction of the concept of interval-like structures. In our search for mathematical contexts where  this kind of structures  naturally arise we have found of some interest to highlight certain kind of them for which the last of the axioms defining an MV algebra (the commutativity of the join) does not hold.

\begin{definition}
A generalized MV-algebra is an MV-algebra for which we do not require the last axiom of MV-algebra.
\end{definition}

We are aware that our definition of generalized conflicts with the one given in \cite{ga}, so we remark the fact that they do not coincide.

Suppose $S$ is a finite set of cardinality $n$;  then we can put a total order on $S$. Any subset of $A\subseteq S$, inherits a total order from $S$, and we indicate its elements, with respect to such an order, by $A_1,A_2,\dots$.

The following example shows that generalized MV algebras can be associated to very basic objects. 

Let $\mathcal{P}(S)$ the power set of $S$: $\mathcal{P}(S)$ is clearly a $\sigma-$ algebra; let $\mu$ the counting measure on $\mathcal{P}(S)$.
\begin{definition}
Given any pair $(A,B)$ in $\mathcal{P}(S)\times\mathcal{P}(S)$ such that $A\subseteq B$, we define 
\[j(A_i)=B\cup(S/B)_{i\ mod|S/B|}\]
\end{definition}
We denote by $j(A)$ the image of entire $A$ under the map $j$.
Then we can define the following operation
\begin{definition}
For any $A,B\ \in\mathcal{P}(S)$, we set 
\[A\oplus B=A\cup B\cup j(A\cap B,A\cup B)\]
\end{definition}
If we indicate with $\neg A$ the complement of $A$, we have the following:
\begin{theorem}
$(\mathcal{P}(S),\oplus,\neg, \emptyset, S)$ is a generalized MV-algebra; for any chain $\mathcal{C}$ of subsets of $S$, $(\mathcal{C},\oplus,\neg,\emptyset,S)$ is an MV algebra.
\end{theorem}

$Proof.$\ The operation is commutative by definition.\\
Given $A,B,C\in \mathcal{S}$ we can check that 
\[(A\oplus B)\oplus C=A\cup B\cup C\cup\Lambda=A\oplus(B\oplus C)\] 
where $\Lambda$ is a subset of $\mathcal{P}(S)$ so defined
\begin{enumerate}
\item $\mu\Lambda=\mu(A\cap B)+\mu(B\cap C)$
\item $\Lambda_i=(S/(A\cup B\cup C))_i$.
\end{enumerate}
It is easy to come up with an example that does not verify the last axiom of MV-algebra. For instance, let us suppose $S=\{a,b,c,d,e\}$, 
A=\{{b,c\} and $B=\{c,d\}$. According to our definition:
\[\neg(\neg A\oplus B)\oplus B=\neg(\{a,d,e\}\oplus\{c,d\})\oplus\{c,d\}=\]
\[= \neg(\{a,d,e,c,d\}\oplus\{c,d\}=\emptyset\oplus\{c,d\}=\{c,d\}=B\]
By the same argument we can deduce that
\[\neg(\neg B\oplus A)\oplus A=A\neq B\]
As we expect, on chain the `join' commutes: let $A\subset B\neq S$; then 
\[\neg(\neg A\oplus B)\oplus B=\neg(S)\oplus B=\emptyset\oplus B=B\]
and 
\[\neg(\neg B\oplus A)\oplus A=\neg(A\cup \neg B)\oplus A =B/A\cup A=B\]
Then we can conclude that 
\[\neg(\neg A\oplus B)\oplus B=\neg(\neg B\oplus A)\oplus A\]

\begin{flushright}
$\square$
\end{flushright}

The class of examples we have proposed suggests that MV- algebras, as generalization of the only one dimensional manifold with disconnected boundary, that is $\mathcal{I}$, arise quite naturally in many contexts, of which the measure theoretical one we have highlighted may well be just one of them. It also suggests that the use of $\mathcal{I}$ and its algebra as mathematical tool to describe topological properties of geometric objects is indeed a promising project of research: in the same way as the unit circle $S^1$ is the archetypal form leading to the definition of fundamental group, MV-algebras could be interpreted as a possible tool to describe some new topological invariant of manifolds.

\end{document}